\documentclass[12pt]{article}
\usepackage{amsfonts}
\usepackage{bbm}
\usepackage{mathrsfs}
\usepackage{amssymb}

\title{An Elimination Lemma for\\
Algebras with PBW Bases}

\author{Huishi Li\thanks{e-mail: huishipp@yahoo.com}\\
{\small Department of Applied Mathematics, College of Information Science and Technology}\\
{\small Hainan University,  Haikou 570228, China}}

\date{}

\begin{document}
\maketitle
\begin{center}
\begin{minipage}{135mm}
{\small {\bf Abstract.} Let $K$ be a field, and $A=K[a_1,\ldots
,a_n]$ a finitely generated $K$-algebra with the PBW $K$-basis
${\cal B}=\{a_{1}^{\alpha_1}\cdots
a_{n}^{\alpha_n}~|~(\alpha_1,\ldots ,\alpha_n)\in\mathbb{N}^n\}$.
It is shown that if $L$ is a nonzero left ideal of $A$ with
GK.dim$(A/L)=d<n$ ($=$ the number of generators of $A$), then $L$
has the {\it elimination property} in the sense that ${\bf V}(U)\cap
L\ne \{0\}$ for every subset $U=\{ a_{i_1},\ldots
,a_{i_{d+1}}\}\subset\{a_1,\ldots ,a_n\}$ with $i_1<i_2<\cdots
<i_{d+1}$, where ${\bf V}(U)=K$-span$\{a_{i_1}^{\alpha_1}\cdots
a_{i_{d+1}}^{\alpha_{d+1}}~|~(\alpha_1,\ldots
,\alpha_{d+1})\in\mathbb{N}^{d+1}\}$. In terms of the structural
properties of $A$, it is also explored when the condition
GK.dim$(A/L)<n$ may hold for a left ideal $L$ of $A$. Moreover, from
the viewpoint of realizing the elimination property by means of
Gr\"obner bases, it is demonstrated that if $A$ is in the class of
binomial skew polynomial rings [G-I2, Serdica Math. J., 30(2004)] or
in the class of solvable polynomial algebras [K-RW, J. Symbolic
Comput., 9(1990)], then every nonzero left ideal $L$ of $A$
satisfies  GK.dim$(A/L)<$ GK.dim$A=n$ ($=$ the number of generators
of $A$), thereby $L$ has the elimination property. }
\end{minipage}\end{center} {\parindent=0pt\vskip 6pt

{\bf MSC 2010} Primary 13P10; Secondary 16W70, 68W30 (16Z05).\vskip
6pt

{\bf Key words} Elimination, PBW $K$-basis, Filtration, 
Gelfand-Kirillov dimension,  Gr\"obner basis.}

\vskip .5truecm

\def\hang{\hangindent\parindent}
\def\textindent#1{\indent\llap{#1\enspace}\ignorespaces}
\def\re{\par\hang\textindent}

\def\v5{\vskip .5truecm}\def\QED{\hfill{$\Box$}}\def\hang{\hangindent\parindent}
\def\textindent#1{\indent\llap{#1\enspace}\ignorespaces}
\def\item{\par\hang\textindent}
\def \r{\rightarrow}\def\OV#1{\overline {#1}}
\def\normalbaselines{\baselineskip 24pt\lineskip 4pt\lineskiplimit 4pt}
\def\mapdown#1{\llap{$\vcenter {\hbox {$\scriptstyle #1$}}$}
                                \Bigg\downarrow}
\def\mapdownr#1{\Bigg\downarrow\rlap{$\vcenter{\hbox
                                    {$\scriptstyle #1$}}$}}
\def\mapright#1#2{\smash{\mathop{\longrightarrow}\limits^{#1}_{#2}}}
\def\NZ{\mathbb{N}}\def\mapleft#1#2{\smash{\mathop{\longleftarrow}\limits^{#1}_{#2}}}

\def\LH{{\bf LH}}\def\LM{{\bf LM}}\def\LT{{\bf
LT}}\def\KX{K\langle X\rangle} \def\KS{K\langle X\rangle}
\def\B{{\cal B}} \def\LC{{\bf LC}} \def\G{{\cal G}} \def\FRAC#1#2{\displaystyle{\frac{#1}{#2}}}
\def\SUM^#1_#2{\displaystyle{\sum^{#1}_{#2}}} \def\O{{\cal O}}  \def\J{{\bf J}}\def\BE{\B (e)}
\def\PRCVE{\prec_{\varepsilon\hbox{-}gr}}\def\BV{\B (\varepsilon )}\def\PRCEGR{\prec_{e\hbox{\rm -}gr}}

\def\KS{K\langle X\rangle}
\def\LR{\langle X\rangle}\def\T#1{\widetilde #1}
\def\HL{{\rm LH}}\def\NB{\mathbb{N}}\def\HY{\hbox{\hskip .03truecm -\hskip .03truecm}}\def\F{{\cal F}}


\section*{0. Introduction}
Let $R=K[x_1,\ldots ,x_n]$ be the commutative polynomial $K$-algebra
in $n$ variables over a field $K$, and let $\{ x_{i_1},\ldots
,x_{i_s}\}\subset\{ x_1,\ldots ,x_n\}$ with $i_1<i_2<\ldots <i_s$.
We say that a monomial ordering on $R$ is an {\it elimination
ordering of type $s$} with respect to the subalgebra
$K[x_{i_1},\ldots ,x_{i_s}]$, denoted $\prec_s$, if $f\in R$ with
the leading monomial $\LM_{\prec_s}(f)\in K[x_{i_1},\ldots
,x_{i_s}]$ implies $f\in K[x_{i_1},\ldots ,x_{i_s}]$.  Let $I$ be a
nonzero ideal of $R$.  Then it follows from Buchberger's Gr\"obner
basis theory that there is the {\it Elimination Theorem for ideals
of $R$}, which states that {\parindent=.44truecm\vskip 6pt

\item{$\bullet$}  If $\G$ is a Gr\"obner basis of $I$ with respect to $\prec_s$,
then $\G_s =\G \cap K[x_{i_1},\ldots ,x_{i_s}]$ is a Gr\"obner basis
of the ideal $I\cap K[x_{i_1},\ldots ,x_{i_s}]$ in $K[x_{i_1},\ldots
,x_{i_s}]$. \vskip 6pt} {\parindent=0pt

Obviously, from this theorem we may see that $\G_s=\emptyset$
$\Leftrightarrow$ $I\cap K[x_{i_1},\ldots ,x_s]=\{0\}.$ So, without
involving  Gr\"obner basis theory, for an arbitrarily given proper
ideal $I$, it is natural to ask {\parindent=.44truecm\vskip 6pt

\item{$\bullet$} To what extent can the elimination of certain variables happen in $I$ via pure structural
properties of an ideal?\vskip 6pt}

As the literature shows, so far perhaps the best answer to the above
question is the one coming  from the dimension theory in commutative
algebraic geometry. More precisely, recall from [Gr\"o, 1968, 1970]
that a subset $U=\{ x_{i_1},\ldots ,x_{i_r}\}\subset \{ x_1,\ldots
,x_n\}$ with $i_1<i_2<\ldots <i_r$ is said to be {\it independent}
(mod $I$) if $I\cap K[x_{i_1},\ldots ,x_{i_r}]=\{0\}$; otherwise $U$
is called {\it dependent} (mod $I$). Considering the dimension
dim${\cal V}(I)$ of the affine algebraic set ${\cal V}(I)$, it is
now well known from the literature (e.g. [KW], [BW]) that
$$\begin{array}{rcl}\hbox{dim}{\cal V}(I)
&=&\max\left\{|U|~\Big |~U\subset\{x_1,\ldots
,x_n\}~\hbox{independent}~(\hbox{mod}~I)\right\}\\
 &=&\hbox{degree of the Hilbert polynomial
of}~R/I.\end{array}$$
Clearly, the above result tells us that \vskip
6pt{\parindent=.44truecm

\re{$\bullet$} if dim${\cal V}(I)=d<n$, then
$K[x_{i_1},...,x_{i_{d+1}}]\cap I\ne \{ 0\}$ for every subset
$U_{d+1}=\{ x_{i_1},...,x_{i_{d+1}}\}\subset$ $\{ x_1,...,x_n\}$
with $i_1<i_2<\cdots <i_{d+1}$, i.e., there are nonzero elements of
$I$ that only depend on the generators in $U_{d+1}$, in particular,
$K[x_1,...,x_d,x_{d+i}]\cap I\ne \{ 0\},~i=1,...,n-d.$\vskip 6pt}
{\parindent=0pt

At this stage, if dim${\cal V}(I)=d<n$ and we want to take out a
nonzero polynomial from $K[x_{i_1},...,x_{i_{d+1}}]\cap I$,  then,
with an elimination ordering $\prec_{d+1}$ of type $d+1$ with
respect to $K[x_{i_1},...,x_{i_{d+1}}]$, running Buchberger's
algorithm will produce a Gr\"obner ~basis $\G$ for $I$ such that
$\G$ contains a nonzero polynomial of $K[x_{i_1},...,x_{i_{d+1}}]$.
}}\par

Here let us point out that in algorithmically computing dim${\cal
V}(I)$  by using a Gr\"obner basis of $I$ ([KW], [BW]) , a {\it
strong  independence} of $U$ (mod $I$) for  a subset $U\subset \{
x_1,\ldots ,x_n\}$  was introduced to make the  {\it key} link
between the  independence of $U$ (mod $I$) and a Gr\"obner basis of
$I$,  and  a {\it a graded monomial ordering that respects the total
degree of polynomials} was necessarily employed throughout the whole
implementing process.\v5

Turning to the noncommutative case, let
$A_n(\mathbb{C})=\mathbb{C}[x_1,\ldots ,x_n,\partial_1,\ldots
\partial_n]$ be the $n$-th Weyl algebra over the field $\mathbb{C}$ of
complex numbers, where $x_1,\ldots ,x_n$ are indeterminate over
$\mathbb{C}$, $\partial_i=\frac{\partial~}{\partial x_i}$, $1\le
i\le n$, and let $L$ be a left ideal of $A_n(\mathbb{C})$. Then the
well-known {\it elimination lemma for  Weyl algebras} [Zei1, Lemma
4.1] states that {\parindent=.44truecm\vskip 6pt

\item{$\bullet$} If
$A_n(\mathbb{C})/L$ is a holonomic $A_n(\mathbb{C})$-module (i.e.
the Gelfand-Kirillov dimension GK.dim$A_n(\mathbb{C})/L=n$), then,
for every $n+1$ generators out of the $2n$ generators $\{
x_1,...,x_n,\partial_1,...,\partial_n\}$ of $A_n(\mathbb{C})$ there
is a nonzero member of $L$ that only depends on these $n+1$
generators. In particular, for every $i=1,...,n$, $L$ contains a
nonzero element of the subalgebra $\mathbb{C}[x_1,...,
x_n,\partial_i]\subset A_n(\mathbb{C})$. \vskip 6pt}{\parindent=0pt

As $A_n(\mathbb{C})$ is a solvable polynomial algebra in the sense
of [KR-W] and it admits the pure lexicographic ordering
$x_1\prec_{lex}\cdots
\prec_{lex}x_n\prec_{lex}\partial_{i_1}\prec_{lex}\cdots\prec_{lex}\partial_{i_n}$
which is certainly an elimination ordering of type $n+1$ with
respect to the subalgebra $\mathbb{C}[x_1,..., x_n,\partial_{i_1}]$,
i.e., $f\in A_n(\mathbb{C})$ with the leading monomial $\LM (f)\in
\mathbb{C}[x_1,..., x_n,\partial_{i_1}]$  implies $f\in
\mathbb{C}[x_1,..., x_n,\partial_{i_1}]$, it follows that if
GK.dim$A_n(\mathbb{C})/L=n$,  then running a noncommutative version
of Buchberger's algorithm constructed in [K-RW] will produce a left
Gr\"obner basis $\G$ for $L$ such that $\G$ contains a nonzero
element of the subalgebra $\mathbb{C}[x_1,..., x_n,\partial_{i_1}]$.
While concerning the determination of holonomicity of
$A_n(\mathbb{C})/L$ ( i.e. GK.dim$A_n(\mathbb{C})/L=n$) by using
Gr\"obner bases, that may refer to a much more general story about
computation of Gelfand-Kirillov dimension for modules over quadric
solvable polynomial algebras [Li1, CH.V] (we shall soon come to this
point below). Nowadays, the elimination lemma for Weyl algebras
[Zei1, Lemma 4.1] has been referred to as the ``fundamental lemma''
in the automatic proving of holonomic function identities [WZ].
Based on this lemma, effective automatic proving  of holonomic
function identities has been carried out, and a large class of
special function identities has been identified ([PWZ], [Zei2],
[Ch], [CS]).}\par

Furthermore, consider a (noncommutative) quadric solvable polynomial
algebra $A=K[a_1,\ldots ,a_n]$ in the sense of [Li1, CH.III, Section
2], that admits a {\it graded monomial ordering $\prec_{gr}$
respecting every $a_i$ being of degree 1}  (Weyl algebras are
typical examples of such algebras). Since  $A$  has the PBW
$K$-basis $\B =\{ a_{1}^{\alpha_1}\cdots
a_{n}^{\alpha_n}~|~(\alpha_1,\ldots ,\alpha_n)\in\NZ^n\}$ (see next
section for an interpretation of this notion), for a subset
$U=\{a_{i_1},\ldots ,a_{i_r}\}\subset\{a_1,\ldots ,a_n\}$ with
$i_1<i_2<\cdots <i_r$,   $U$ is said to be {\it weakly independent
modulo a left ideal  $L$ of $A$} if $L\cap {\bf V}(U)=\{0\}$, where
$${\bf V}(U)=K\hbox{-span}\left\{\left. a_{i_1}^{\alpha_1}\cdots
a_{i_r}^{\alpha_r}~\right |~(\alpha_1,...,\alpha_r)
\in\NZ^r\right\}.$$
With this weak independence of $U$ (mod $L$) and   a double
filtered-graded transfer trick,  the strategy of computing dim${\cal
V}(I)$ proposed by ([KW], [BW]) was adapted in [Li1, CH.V] for
computing the Gelfand-Kirillov dimension GK.dim$(A/L)$, and
consequently the following results were established: \vskip 6pt
{\parindent=.44truecm

\item{$\bullet$}  [Li1, CH.V, Theorem  7.4]  Let $L$ be a {\it nonzero} left
ideal of $A$. Then
$$\begin{array}{rcl}\hbox{GK.dim}(A/L)
&=&\hbox{degree of the Hilbert polynomial of}~A/L\\
&=&\max\left\{|U|~\Big |~U\subset\{a_1,\ldots ,a_n\}~\hbox{weakly
independent}~(\hbox{mod}~L)\right\},\end{array}$$  which can be
algorithmically computed  via a Gr\"obner basis of $L$; moreover,
$$\hbox{GK.dim}(A/L)< n=\hbox{GK.dim}A;$$

\item{$\bullet$} [Li1, CH.V, Lemma 7.5]  If GK.dim$A/L=d$, then ${\bf V}(U)\cap L\ne \{ 0\}$
for every subset $U=\{ a_{i_1},...,a_{i_{d+1}}\}\subset \{
a_1,...,a_n\}$ with $i_1<i_2<\cdots <i_{d+1}$.\vskip
6pt}{\parindent=0pt Therefore, if GK.dim$A/L=d$, $U=\{
a_{i_1},...,a_{i_{d+1}}\}$ with $i_1<i_2<\cdots <i_{d+1}$, and if
$A$ admits an elimination ordering $\prec_{d+1}$ of type $d+1$ with
respect to ${\bf V}(U)$, i.e., $f\in A$ with the leading monomial
$\LM (f)\in {\bf V}(U)$ implies $f\in {\bf V}(U)$, then running the
noncommutative Buchberger's algorithm constructed [K-RW] will
produce a left Gr\"obner ~basis $\G$ for $L$ such that $\G$ contains
a nonzero element of ${\bf V}(U)$.\par}

Note that the class of quadric solvable polynomial algebras studied
in [Li1, CH.III, CH.V] covers not only Weyl algebras, but also more
(skew) Ore extensions and operator algebras. Enlightened by the
automatic proving of multivariate identities over operator algebras
([PWZ], [Ch], [CS]),   more general  $\partial$-finiteness and
$\partial$-holonomicity for modules over quadric solvable polynomial
algebras were introduced and preliminarily studied  in [Li1. CH.VII]
by using [Li1, CH.V, Lemma 7.5] as a key role. \v5

Inspired by the elimination lemma [Zei1, Lemma 4.1] and the
elimination Lemma [Li1, CH.V, Lemma 7.5], in this paper we first
show that there is a kind of elimination lemma (Elimination Lemma
2.1) for {\it all} finitely generated $K$-algebras with PBW
$K$-bases, that is,  if $A=K[a_1,\ldots ,a_n]$ is a finitely
generated $K$-algebra with the PBW $K$-basis $\B =\{
a_{1}^{\alpha_1}\cdots a_{n}^{\alpha_n}~|~(\alpha_1,\ldots
,\alpha_n)\in\NZ^n\}$,  and if $L$ is a nonzero left ideal of $A$
with GK.dim$(A/L)=d<n$ ($=$ the number of generators of $A$), then
$L$ has the {\it elimination property} in the sense that ${\bf
V}(U)\cap L\ne \{0\}$ for every subset $U=\{ a_{i_1},\ldots
,a_{i_{d+1}}\}\subset\{a_1,\ldots ,a_n\}$ with $i_1<i_2<\cdots
<i_{d+1}$, where ${\bf V}(U)=K$-span$\{a_{i_1}^{\alpha_1}\cdots
a_{i_{d+1}}^{\alpha_{d+1}}~ |~(\alpha_1,\ldots
,\alpha_{d+1})\in\mathbb{N}^{d+1}\}$; then we explore, in terms of
the structural properties of $A$, when the condition GK.dim$(A/L)<n$
may hold for a nonzero left ideal $L$ of $A$ (Theorem 2.6). From the
viewpoint of realizing the elimination property by means of
Gr\"obner bases, in the last two sections, we demonstrate that if
$A=K[a_1,\ldots ,a_n]$ is an algebra in either the class of binomial
skew polynomial rings [G-I2] or the class of solvable polynomial
algebras [K-RW], then every nonzero left ideal $L$ of $A$ satisfies
GK.dim$(A/L)<$ GK.dim$A=n$, thereby $L$ has the elimination property
(Theorem 3.3, Theorem 4.1). \v5

Throughout the following sections, $K$ denotes a field, $K^*=K-\{
0\}$; $\mathbb{N}$ denotes the set of all nonnegative integers, and
$\mathbb{Z}$ denotes the set of all integers; algebras are meant
associative $K$-algebra with multiplicative identity 1; if
$A=K[a_1,\ldots ,a_n]$ is a finitely generated $K$-algebra, then we
always assume that the set of generators $\{ a_1,\ldots ,a_n\}$ is
{\it minimal}, i.e.,  any proper subset of  $\{ a_1,\ldots ,a_n\}$
cannot generate $A$ as a $K$-algebra.

\section*{1. Preliminaries}
To reach our goal of this paper, in this section we recall several
necessary  notions  concerning  finitely generated $K$-algebras and
their modules, such as PBW $K$-basis, $\NZ$-filtration, and
Gelfand-Kirillov dimension; moreover, some known results related to
these notions, which will be used in later sections, are recalled as
well.  A general Gr\"obner basis theory for ideals of free algebras
is referred to [Gr] and [Mor], and a general Gelfand-Kirillov
dimension theory for algebras and modules is referred to [KL] and
[MR].\v5

Let $A=K[a_1,\ldots ,a_n]$ be a finitely generated $K$-algebra with
the  set of generators $\{ a_1,\ldots ,a_n\}$. If, for some
permutation $i_1i_2\cdots i_n$ of $1,2,\ldots ,n$, the set $\B =\{
a^{\alpha}=a_{i_1}^{\alpha_1}\cdots a_{i_n}^{\alpha_n}~|~\alpha
=(\alpha_1,\ldots ,\alpha_n)\in\NZ^n\} ,$ forms a $K$-basis of $A$,
then $\B$ is referred to as a {\it PBW $K$-basis} of $A$ (where the
phrase ``PBW $K$-basis" is  abbreviated from the well-known {\it
Poincar\'e-Birkhoff-Witt Theorem} concerning the standard $K$-basis
of the enveloping algebra of a Lie algebra, e.g. see [Hu, P. 92]).
For more  content concerning PBW $K$-bases and related topics, the
reader is referred to a nice survey paper [SW].\par

It is clear that if $A$ has a PBW $K$-basis, then we can always
assume that $i_1=1,\ldots ,i_n=n$. Thus, we make the following
convention once for all.{\parindent=0pt\v5

{\bf Convention} From now on in this paper, if we say that a
finitely generated algebra $A=K[a_1,\ldots ,a_n]$ has the PBW
$K$-basis $\B$, then it always means that
$$\B =\{ a^{\alpha}=a_{1}^{\alpha_1}\cdots
a_{n}^{\alpha_n}~|~\alpha =(\alpha_1,\ldots ,\alpha_n)\in\NZ^n\} .$$
Moreover, adopting the commonly used terminology in computational
algebra, elements of $\B$ are referred to as {\it monomials} of
$A$.}\v5

Let $\KX =K\langle X_1,\ldots ,X_n\rangle$ be the free $K$-algebra
generated by $X=\{ X_1,\ldots ,X_n\}$. Then $\KX$ has the standard
$K$-basis $\mathbb{B}$ consisting of words on alphabet $X$, or more
precisely, writing 1 for the empty word,
$$\mathbb{B}=\{1\}\cup\{X_{i_1}\cdots X_{i_s}~|~X_{i_j}\in X,~s\ge 1\} .$$
Since every finitely generated $K$-algebra $A=K[a_1,\ldots ,a_n]$
has a presentation $\KX /I$ with respect to $X_i\mapsto a_i$, $1\le
i\le n$, where $I$ is an ideal of $\KX$, the proposition stated
below, which  is a generalization of [Gr, Proposition 2,14] and
[Li1, CH.III, Theorem 1.5], may be viewed as an algorithmic
criterion for $A$ to have the PBW $K$-basis $\B$.{\parindent=0pt\v5

{\bf 1.1. Proposition} [Li4, Ch 4, Theorem 3.1] Let $I\ne \{ 0\}$ be
an ideal of the free $K$-algebra $\KX =K\langle X_1,\ldots
,X_n\rangle$, and $A=\KS /I$. Suppose that $I$ contains a subset of
$\frac{n(n-1)}{2}$ elements
$$G=\{ g_{ji}=X_jX_i-F_{ji}~|~F_{ji}\in\KS ,~1\le i<j\le n\}$$
such that with respect to some monomial ordering $\prec_{_X}$ on
$\mathbb{B}$, the leading monomial $\LM (g_{ji})=X_jX_i$ for all
$g_{ji}\in G$. The following two statements are equivalent:\par

(i) $A$ has the  PBW $K$-basis $\B =\{ \OV X_1^{\alpha_1}\OV
X_2^{\alpha_2}\cdots \OV X_n^{\alpha_n}~|~\alpha_j\in\NZ\}$ where
each $\OV X_i$ denotes the coset of $I$ represented by $X_i$ in
$A$.\par

(ii) Any  subset $\G$ of $I$ containing $G$ is a  Gr\"obner basis
for $I$ with respect to $\prec_{_X}$.\QED \v5

{\bf Remark}  Obviously, Proposition 1.1 holds true if we use any
permutation $X_{k_1},X_{k_2},\ldots ,X_{k_n}$ of the generators
$X_1,X_2,\ldots ,X_n$ of $\KX$ (see an example given in Section 3).
So, in what follows we conventionally keep using the natural
permutation $X_1,\ldots ,X_n$.}\v5

Note that the multiplication of elements in the $K$-basis
$\mathbb{B}$ of $\KX$ is given by the concatenation of words, i.e.,
$$(X_{i_1}\cdots X_{i_s})\cdot (X_{j_1}\cdots X_{j_t})=X_{i_1}\cdots
X_{i_s}X_{j_1}\cdots X_{j_t}.$$ It follows that $\KX$ is turned into
an $\NZ$-filtered algebra by the natural filtration  $F\KX =\{
F_m\KX\}_{m\in\NZ}$ with
$$
F_m\KX =K\hbox{-span}\{1,~X_{i_1}\cdots X_{i_s}\in\B~|~s\le m\}
,\quad m\in\NZ,$$ which satisfies that $F_0\KX =K$, each $F_m\KX$ is
a finite dimensional subspace of $\KX$,  $\KX
=\cup_{m\in\NZ}F_m\KX$, $F_m\KX\subseteq F_{m+1}\KX$, and
$F_{m_1}\KX F_{m_2}\KX\subseteq F_{m_1+m_2}\KX$. If
$A=K[a_1,...,a_n]$ is a finitely generated $K$-algebra, then there
is an ideal $I$ of $\KX$ and an algebra isomorphism
$\KX/I~\mapright{\cong}{}~A$ with $\OV{X}_i\mapsto a_i$, $1\le i\le
n$, where each $\OV X_i=X_i+I$ is the coset represented by $X_i$ in
$\KX/I$.  Consequently, with respect to $\OV X_i\mapsto a_i$, $1\le
i\le n$, $A$ is turned into an  $\NZ$-filtered algebra by the
filtration $FA=\{ F_mA\}_{m\in\NZ}$ induced by $F\KX$, i.e., for
every $m\in\NZ$,
$$\begin{array}{rcl}
F_mA&=&K\hbox{-span}\left\{1,~a_{i_1}\cdots a_{i_s}~ |~
a_{i_j}\in\{a_1,\ldots ,a_n\},~1\le s\le m\right\} \\
&\cong&(F_m\KX+I)/I,\end{array},$$  which satisfies that $F_0A=K$,
each $F_mA$ is a finite dimensional subspace of $A$,
$A=\cup_{m\in\NZ}F_mA$, $F_mA\subseteq F_{m+1}A$,
$F_{m_1}AF_{m_2}A\subseteq F_{m_1+m_2}A$  for all $m_1,m_2\in\NZ$
(note that this is determined by the concatenation of ``words"
$(a_{i_1}\cdots a_{i_s})\cdot (a_{k_1}\cdots a_{k_t})=a_{i_1}\cdots
a_{i_s}a_{k_1}\cdots a_{k_t}$). Indeed, a direct verification shows
that if we take $V=\sum^n_{i=1}Ka_i$, then $F_mA=K+V+V^2+\cdots
+V^m$ for all $m\in\NZ$. Hence, in the literature (cf. [MR, P.26]),
the $\NZ$-filtration $FA$ as described above is usually referred to
as the {\it standard filtration} of $A$ determined by the {\it
finite dimensional generating subspace} $V$. Considering the
function $d_F(m)=$ dim$_KF_mA$, $m\in\NZ$, the Gelfand-Kirillov
dimension (abbreviated GK dimension) of $A$, denoted GK.dim$A$, is
defined as
$$\hbox{GK.dim}A=\inf\left\{\lambda\in\mathbb{R}~\left |~d_F(m)\le m^{\lambda}~\hbox{for}~m\gg 0\right.\right\}.$$
It is known that GK.dim$A$ does not depend on the choice of a finite
dimensional generating subspace of $A$. So, by the definition,
GK.dim$A$ amounts to a measure of the growth rate of $A$ {\it  as a
$K$-algebra with respect to any finite dimensional generating
subspace}. If the ``$\inf$" exists in the above definition, say
GK.dim$A=\lambda$, then we say that $A$ has polynomial growth; if
the  ``$\inf$" does not exist, then we write GK.dim$A=\infty$.\par

By means of a Gr\"obner basis technique, the next proposition tells
us when a finitely generated $K$-algebra $A=K[a_1,\ldots ,a_n]$ of
$n$ generators has GK.dim$A=n$. {\parindent=0pt\v5

{\bf 1.2. Proposition} [Li2, Section 6, Example 1] (or [Li4, Ch. 5,
Section 5.3, Example 3]) Let the $K$-algebra $A=K[a_1,\ldots ,a_n]$
be presented as a quotient algebra of the free $K$-algebra $\KX
=K\langle X_1,\ldots ,X_n\rangle$, say $A=\KX /I$, and suppose that
the ideal $I$ of $\KX$ has a finite Gr\"obner basis $\G=\{
g_{ji}~|~1\le i<j\le n\}$ with respect to some monomial ordering
$\prec_{_X}$ on the $K$-basis $\mathbb{B}$ of $\KX$, such that the
leading monomial $\LM (g_{ji})=X_jX_i$, $1\le i<j\le n$. Then  $A$
has polynomial growth and GK.dim$A=n=$ the number of generators of
$A$.\par\QED}\v5

Finally we recall the notion of Gelfand-Kirillov dimension (GK
dimension for short) for finitely generated modules over a finitely
generated $K$-algebra $A=K[a_1,\ldots ,a_n]$. Let
$M=\sum^s_{i=1}A\xi_i$ be a  finitely generated left $A$-module with
the generating set $\{\xi_1,\ldots ,\xi_s\}$.  Considering $A$ as an
$\NZ$-filtered algebra with the standard filtration $FA=\{
F_mA\}_{m\in\NZ}$, $M$ is then turned into an $\NZ$-{\it filtered
$A$-module} by the filtration $FM=\{F_qM\}_{q\in\NZ}$ with
$$F_qM=F_qAF_0M,~\hbox{where}~F_0M=\sum^s_{i=1}K\xi_i,~q\in\NZ,$$
which satisfies that each $F_qM$ is a finite dimensional subspace of
$M$, $M=\cup_{q\in\NZ}F_qM$, $F_qM\subseteq F_{q+1}M$,
$F_{m}AF_{q}M\subseteq F_{m+q}M$ for all $m,q\in\NZ$. In the
literature (cf. [MR, P.300]), the filtration $FM$ as described above
is also referred to as the {\it standard filtration} of $M$
determined by the {\it finite dimensional generating subspace}
$\sum^s_{i=1}K\xi_i$. Considering the function $d_F(q)=$
dim$_KF_qM$, $q\in\NZ$, the GK dimension of $M$, denoted GK.dim$M$,
is defined as
$$\hbox{GK.dim}M=\inf\left\{\lambda\in\mathbb{R}~\left |~d_F(q)\le q^{\lambda}~\hbox{for}~q\gg 0\right.\right\}.$$
It is also known that GK.dim$M$ does not depend on the choice of a
finite dimensional generating subspace of $M$. So, by the
definition, GK.dim$M$ amounts to a measure of the growth rate of $M$
{\it as an $A$-module with respect to any finite dimensional
generating subspace}. If the ``$\inf$" exists in the above
definition, say GK.dim$M=\lambda$, then we say that $M$ has
polynomial growth; if the  ``$\inf$" does not exist, then we write
GK.dim$M=\infty$.{\parindent=0pt\v5

\section*{2. Elimination Lemma for Algebras with PBW Bases}
Let $A=K[a_1,\ldots ,a_n]$ be a  finitely generated $K$-algebra with
the PBW $K$-basis $\B =\left\{\left. a^{\alpha}=a_1^{\alpha_1}\cdots
a_n^{\alpha_n}~\right |~\alpha =(\alpha_1,\ldots
,\alpha_n)\in\NZ^n\right\}.$ Consider the positive-degree function
$d$ on $\B$, which assigns $d(a_i)=1$, $1\le i\le n$. Then, for each
$a^{\alpha}\in\B$ with $\alpha =(\alpha_1,\ldots
,\alpha_n)\in\NZ^n$, we have $d(a^{\alpha})=\alpha_1+\cdots
+\alpha_n$. Thus the $K$-{\it vector space} $A$ has the
$\NZ$-filtration $\F A=\{ \F_mA\}_{m\in\NZ}$ with
$$\F_mA=K\hbox{-span}\left\{\left. a^{\alpha}\in\B~\right |~d(a^{\alpha})\le m\right \},\quad m\in\NZ,$$ which satisfies
that $\F_0A=K$,  each $\F_mA$ is a finite dimensional subspace of
$A$, $A=\cup_{m\in\NZ}\F_mA$, and $\F_mA\subseteq \F_{m+1}A$ for all
$m\in\NZ$, but does not necessarily satisfy
$\F_{m_1}A\F_{m_2}A\subseteq \F_{m_1+m_2}A$ for $m_1,m_2\in\NZ$ , or
in other words, $A$ {\it is not necessarily an $\NZ$-filtered
algebra with respect to the filtration} $\F A$ (see the example
given in Section 4). Furthermore, comparing the filtration $\F A=\{
\F_mA\}_{m\in\NZ}$ with the standard filtration
$FA=\{F_mA\}_{m\in\NZ}$ of $A$ defined in the last section,  it is
clear that $\F_mA\subseteq F_mA,\quad m\in\NZ .$}\par

Now, for any subset $U_r=\{ a_{i_1},...,a_{i_r}\}\subset\{
a_1,...,a_n\}$ with $i_1<i_2<\cdots <i_r$,  let us take
$${\bf V}(U_r)=K\hbox{-span}\left\{ a_{i_1}^{\alpha_1}\cdots a_{i_r}^{\alpha_r}~\Big |~
(\alpha_1,...,\alpha_r)\in\NZ^r\right\}.$$ Then we are ready to
establish an Elimination Lemma for the algebra $A$.{\parindent=0pt
\v5

{\bf 2.1. Elimination Lemma} Let $A$ and the notations be as fixed
above, and  let $L$ be a nonzero left ideal of $A$. If the left
$A$-module $A/L$ has finite GK dimension $\hbox{GK.dim}(A/L)=d$,
then,  for any subset  $U_r=\{ a_{i_1},\ldots ,a_{i_r}\}
\subset\{a_1,\ldots ,a_n\}$ with $i_1<i_2<\cdots <i_r$,   ${\bf
V}(U_r)\cap L=\{ 0\}$ implies $r\le d$.  Consequently, if $d<n$ ($=$
the number of generators of $A$), then  ${\bf V}(U_{d+1})\cap L\ne
\{ 0\}$ holds true for every subset $U_{d+1}=\{
a_{i_1},...,a_{i_{d+1}}\}\subset$ $\{ a_1,...,a_n\}$ with
$i_1<i_2<\cdots <i_{d+1}$, in particular, for every $U_s=\{
a_1,\ldots a_s\}$ with $d+1\le s\le n-1$ we have ${\bf V}(U_s)\cap
L\ne \{ 0\}$. \vskip 6pt

{\bf Proof}  Let the $A$-module $A/L$ be equipped with the
filtration $F(A/L)=\{ F_q(A/L)\}_{q\in\NZ}$ induced by the standard
filtration $FA=\{F_qA\}_{q\in\NZ}$ of $A$, i.e.,
$F_q(A/L)=(F_qA+L)/L$, $q\in\NZ$, which is clearly the standard
filtration of $A/L$ as defined in Section 1. Taking any $U_r=\{
a_{i_1},...,a_{i_r}\}\subset\{ a_1,...,a_n\}$ with $i_1<i_2<\cdots
<i_r$, consider the filtration $\F {\bf V}(U_r)=\{\F_q{\bf
V}(U_r)\}_{q\in\NZ}$ of the vector space ${\bf V}(U_r)$ induced by
the filtration $\F A=\{\F_qA\}_{q\in\NZ}$ of $A$ determined by the
PBW basis $\B$ (as defined above), i.e., $\F_q{\bf V}(U_r)={\bf
V}(U_r)\cap\F_qA$, $q\in\NZ$. If ${\bf V}(U_r)\cap L=\{0\}$, then
since
$$\F_q{\bf V}(U_r)\cong \frac{\F_q{\bf V}(U_r)+L}{L}\subset\frac{\F_qA+L}{L}\subset \frac{F_qA+L}{L},\quad q\in\NZ,$$
it turns out  that for every $q\in\NZ$,
$$\begin{array}{rcl} \left (\begin{array}{c} q+r\\ r\end{array}\right )&=&\hbox{dim}_K\F_q{\bf V}(U_r)\\
&=&\hbox{dim}_K\displaystyle{\frac{\F_q{\bf V}(U_r)+L}{L}}\\
\\
&\le&\hbox{dim}_K\displaystyle{\frac{F_qA+L}{L}}\\
\\
&=&\hbox{dim}_KF_q(A/L).\end{array}$$ Now if $\hbox{GK.dim}M=d$,
then, taking the usual infimum
$$\inf\{\lambda\in\mathbb{R}~|~\hbox{dim}_K\F_q{\bf V}(U_r)\le
q^{\lambda},~q\gg 0\}$$ into account (of course this infimum is
nothing about GK dimension), it follows from  the definition of GK
dimension for the module $M$ that
$$\begin{array}{rcl} r&=&\inf\{\lambda\in\mathbb{R}~|~\hbox{dim}_K\F_q{\bf V}(U_r)\le q^{\lambda},~q\gg 0\}\\
&\le&\inf\{\lambda\in\mathbb{R}~|~\hbox{dim}_KF_q(A/L)\le
q^{\lambda},~q\gg 0\}\\
&=&\hbox{GK.dim}(A/L)=d,\end{array}$$ as desired. Consequently, the
last assertion of the lemma is immediately clear.\QED}\v5

For convenience in using Elimination Lemma 2.1, it is reasonable to
introduce  {\parindent=0pt\v5

{\bf 2.2. Definition} Let $A=K[a_1,\ldots ,a_n]$ be a finitely
generated $K$-algebra with the PBW $K$-basis, and let $L$ be a
nonzero left ideal of $A$. If
$$(*)\quad \hbox{GK.dim}(A/L)<n$$
then we say that {\it Elimination Lemma 2.1 holds true for $L$}.
}\v5

Since by [Li1, CH.V, Theorem 7.4] we know that every nonzero left
ideal $L$ of a quadric solvable polynomial algebra $A=K[a_1,\ldots
,a_n]$ satisfies $\hbox{GK.dim}(A/L)<$ GK.dim$A=n$, it follows that
$L$ satisfies the condition $(*)$ of Definition 2.2, and hence
Elimination Lemma 2.1 holds true for every nonzero left ideal $L$ of
$A$. Thereby Elimination Lemma 2.1 covers the elimination lemma for
quadric solvable polynomial algebras [Li1, Lemma 7.5], especially it
covers  the elimination lemma for Weyl algebras [Zei, Lemma 4.1]
(note that Weyl algebras are typical quadric solvable polynomial
algebras).\v5

Let $A$ be an arbitrary finitely generated $K$-algebra with the PBW
$K$-basis $\B$, and $L$ a nonzero left ideal of $A$. From a
practical viewpoint, it seems that computing GK.dim$(A/L)$ in an
algorithmic way (as in [Li1, CH.V]) is not always feasible in order
to  realize the condition $(*)$ of Definition 2.2 for $L$. Instead,
learning from the knowledge of Gelfand-Kirillov dimension for
algebras and modules, we would rather try to  determine whether
$$\hbox{GK.dim}(A/L)<\hbox{GK.dim}A\le  n ~(=\hbox{the number of generators of}~A).$$\par

As to getting the first inequality
$\hbox{GK.dim}(A/L)<\hbox{GK.dim}A$,  the  lemma  presented below
may  shed light on this topic. {\parindent=0pt\v5

{\bf 2.3. Lemma} Let $A$ be any  $K$-algebra. Then the following
statements hold.\par

(i) If $f\in A$ is not a divisor of zero, then  GK.dim$(A/Af) \le$
GK.dim$A-1$.\par

(ii) If $A$ is a domain and $L$ is any nonzero left ideal of $A$,
then GK.dim$(A/L)\le $ GK.dim$A-1$.\vskip 6pt

{\bf Proof} (i) If $f\in A$ is not a divisor of zero, then $A\cong
Af$ as left $A$-modules. It follows from [MR, Proposition 8.3.5]
that GK.dim$(A/Af) \le $ GK.dim$A -1$.\par

(ii) If $A$ is a domain and $L$ is any nonzero left ideal of $A$,
then, taking a nonzero $f\in L$, It follows from (i) and the exact
sequence $A/Af\r A/L\r 0$ of $A$-modules that GK.dim$(A/L)\le$
GK.dim$(A/Af)\le$ GK.dim$A-1$. \QED} \v5

Concerning the second inequality $\hbox{GK.dim}A\le  n ~(=\hbox{the
number of generators of}~A),$ it is certainly a matter of
determining GK.dim$A$. As one may know from the literature, there
are many different ways to determine the GK dimension of an algebra.
Especially for a finitely generated $K$-algebra $A=K[a_1,\ldots
,a_n]$, if $A$ is presented as a quotient algebra $\KX /I$ of the
free $K$-algebra $\KX =K\langle X_1,\ldots ,X_n\rangle$ and if the
ideal $I$ has a finite Gr\"obner basis $\G$ with respect to some
monomial ordering $\prec$, then it follows from [Uf] that GK.dim$A$
can be read out from the Ufnarovski  graph of $\G$ (see also [Li2,
Section 6] and [Li4, Ch.5] for a number of examples including the
foregoing Proposition 1.2).  While for an arbitrary finitely
generated $K$-algebra $A$ with the PBW $K$-basis, the next lemma may
also help us to determine the GK dimension of $A$.
{\parindent=0pt\v5

 {\bf 2.4. Lemma} Let $A=K[a_1,\ldots ,a_n]$ be a
finitely generated $K$-algebra with the PBW $K$-basis $\B$, and let
$\F A=\{ \F_mA\}_{m\in\NZ}$ be the $\NZ$-filtration of the vector
space $A$ determined by $\B$ as described before, i.e.,
$\F_mA=K\hbox{-span}\left\{\left. a^{\alpha}\in\B~\right
|~d(a^{\alpha})\right.$ $\left.\le m\right \},$ $m\in\NZ$. Then, the
following two statements are equivalent:\par

(i) The generators of $A$ satisfy
$$a_ja_i =\sum_{q\le
\ell}\lambda_{q\ell}a_qa_{\ell}+\sum_t\lambda_ta_t+\lambda_{ji},~
1\le i<j\le n,~\lambda_{q\ell},\lambda_t,\lambda_{ji}\in K.$$

(ii)  $A$ is turned into an $\NZ$-filtered algebra by $\F A$, i.e.,
$\F_{m_1}A\F_{m_2}A\subseteq \F_{m_1+m_2}A$ for all
$m_1,m_2\in\NZ$,\vskip 6pt

{\bf Proof} Note that the filtration  $\F A$ is constructed by using
the  positive-degree function $d(~)$ on $\B$ such that  $d(a_i)=1$,
$1\le i\le n$. So it is straightforward to verify that $A$ is turned
into an $\NZ$-filtered algebra by $\F A$ if and only if  $a_ja_i$
has the desired representation.\par\QED\v5

{\bf 2.5. Proposition} Let $A=K[a_1,\ldots ,a_n]$ be a finitely
generated $K$-algebra with the PBW $K$-basis $\B$. Consider the
standard filtration  $FA=\{ F_mA\}_{m\in\NZ}$ of the $\NZ$-filtered
algebra $A$, and the $\NZ$-filtration $\F A=\{ \F_mA\}_{m\in\NZ}$ of
the vector space $A$ determined by the PBW $K$-basis  $\B$, as
described before.  If $A$ is also an $\NZ$-filtered algebra with
respect to the filtration $\F A$, then $F_mA=\F_mA$ for all
$m\in\NZ$, and consequently GK.dim$A=n$. \vskip 6pt

{\bf Proof}  By the construction of both $FA$ and $\F A$, it is
clear that $\F_mA\subseteq F_mA$ holds for all $m\in\NZ$. On the
other hand, if $A$ is also an $\NZ$-filtered algebra with respect to
the filtration $\F A$, then  applying Lemma 2.4 (i) to the structure
of both filtration $FA$ and $\F A$, the inclusion
$F_mA\subseteq\F_mA$ is obtained  for every $m\in\NZ$. It follows
that $F_mA=\F_mA$ for all $m\in\NZ$, and consequently
$$\left (\begin{array}{c} m+n\\ n\end{array}\right ) =\hbox{dim}_K\F_mA=\hbox{dim}_KF_mA,\quad m\in\NZ.$$
Thereby the knowledge of Gelfand-Kirillov dimension for algebras
entails
$$\begin{array}{rcl} n&=&\inf\{\lambda\in\mathbb{R}~|~\hbox{dim}_K\F_mA\le m^{\lambda},~m\gg 0\}\\
&=&\inf\{\lambda\in\mathbb{R}~|~\hbox{dim}_KF_mA\le
m^{\lambda},~m\gg 0\}\\
&=&\hbox{GK.dim}A,\end{array}$$ as desired. \QED}\v5

With the aid of the foregoing preparation, the next theorem provides
us with a class of algebras such that if an algebra $A$ belongs to
this class, then every nonzero left ideal of $A$ satisfies the
condition $(*)$ of Definition 2.2.{\parindent=0pt\v5

{\bf 2.6. Theorem} Let $A=K[a_1,\ldots ,a_n]$ be a finitely
generated $K$-algebra. If {\parindent=1.6truecm\par

\item{(1)} $A$ has the PBW $K$-basis $\B$,\par

\item{(2)} the generators of $A$ satisfy
$$a_ja_i =\sum_{q\le
\ell}\lambda_{q\ell}a_qa_{\ell}+\sum_t\lambda_ta_t+\lambda_{ji},~
1\le i<j\le n,~\lambda_{q\ell},\lambda_t,\lambda_{ji}\in K,$$\par}
and {\parindent=1.6truecm

\item{(3)} $A$ is a domain, \par}

then GK.dim$A=n$, and Elimination lemma 2.1 holds true (in the sense
of Definition 2.2) for every nonzero left ideal $L$ of $A$. \vskip
6pt

{\bf Proof} This is just a result of combining Lemma 2.3, Lemma 2.4
and Proposition 2.5.\QED \v5

{\bf Remark} In the next two sections,  we shall respectively
determine two significant subclasses of the class of finitely
generated $K$-algebras with the PBW $K$-basis $\B$, such that if an
algebra $A$ belongs to either of the two subclasses then Elimination
Lemma holds true (in the sense of Definition 2.2) for every nonzero
left ideal $L$ of $A$; moreover, the two subclasses of algebras
will also illustrate that

(i) if a finitely generated $K$-algebra $A=K[a_1,\ldots ,a_n]$ is
presented as a quotient algebra of the free $K$-algebra $\KX
=K\langle X_1,\ldots ,X_n\rangle$, i.e., $A=\KX /I$, then to a large
extent, the Gr\"obner basis technique as shown in Proposition 1.1
will be quite helpful for us to check whether the conditions (1) and
(2) of Theorem 2.6 are satisfied by $A$;
\par

(ii) the class of algebras satisfying the three conditions of
Theorem 2.6  properly contains the class of all quadric solvable
polynomial algebras in the sense of [Li1, CH.III, Section 2];\par

(iii) a finitely generated $K$-algebra  $A$ satisfying the
conditions (1) and (3) of Theorem  2.6 and such that  Elimination
lemma 2.1 holds true (in the sense of Definition 2.2) for every
nonzero left ideal,  may not necessarily  satisfy the condition (2)
of Theorem 2.6 (see the example given in Section 4). }\v5

\section*{3. An Application to Binomial Skew Polynomial Rings}
In the algebraic study of finding solutions to Yang-Baxter
equations, the class of  binomial skew polynomial rings was
introduced and studied in [G-I1,2], and quite rich results were
obtained in which the most important results are: every binomial
skew-polynomial ring A is respectively {\parindent=1.3truecm\par

\item{(1)} a left and right Noetherian domain;\par

\item{(2)} an Artin-Schelter regular PBW algebra;\par

\item{(3)} a Koszul algebra such that the Koszul dual $A!$ is a quantum
Grassmann algebra;\par

\item{(4)} a quantum binomial PBW algebra in the sense of [G-I3], and
hence a Yang-Baxter algebra, that is, the set of defining relations
${\cal R}$ of $A$ defines canonically a solution to the Yang-Baxter
equation.\par}{\parindent=0pt

In this section we demonstrate that every  binomial skew-polynomial
ring $A$ satisfies the three conditions of Theorem 2.6. To better
understand this result, with notions and notations as used in
previous sections, we start by recalling from loc cit. the
definition of a binomial skew-polynomial ring.\v5

{\bf 3.1. Definition} Let $\KX =\langle X_1,\ldots ,X_n\rangle$ be
the free $K$-algebra generated by $X=\{ X_1,\ldots ,X_n\}$ with the
standard $K$-basis $\mathbb{B}=\{ 1\}\cup\{X_{i_1}\cdots
X_{i_s}~|~X_{i_j}\in X,~s\ge 1\}$, and let $I$ be an ideal of $\KX$,
$A=\KX /I$. If every $X_i$ is assigned the degree 1, and if $I$ is
generated by the subset ${\cal R}=\{
R_{ji}=X_jX_i-c_{ij}X_{i'}X_{j'}\}_{1\le i<j\le n}$ consisting of
exactly $\frac{n(n-1)}{2}$ elements, such
that{\parindent=1.3truecm\par

\item{(a)} $c_{ij}\in K^*=K-\{0\}$, $1\le i<j\le n$;\par

\item{(b)} each $R_{ji}=X_jX_i-c_{ij}X_{i'}X_{j'}$ satisfies $i'<j$, $i'< j'$, $1\le i<j\le n$;\par

\item{(c)} $\{ X_{i'}X_{j'}~|~X_jX_i-c_{ij}X_{i'}X_{j'}=R_{ji}\in{\cal R}\}=\{X_iX_j~|~1\le i<j\le n\}$;\par

\item{(d)} with respect to the graded lexicographic monomial ordering $X_1\prec_{grlex}X_2\prec_{grlex}\cdots\prec_{grlex}
X_n$, ${\cal R}$ forms a reduced Gr\"obner basis of $I$.\par}

then $A$ is called a {\it  binomial skew polynomial ring}. \v5

{\bf 3.2. Theorem} Let $A=\KX /I$ be a binomial skew polynomial ring
as defined above. Then  $A$ satisfies the three conditions of
Theorem 2.6, thereby  Elimination lemma 2.1 holds true (in the sense
of Definition 2.2) for every nonzero left ideal $L$ of $A$. \vskip
6pt

{\bf Proof} Though this assertion may follow directly from the
definition and the structural properties of a binomial skew
polynomial ring as we listed in the beginning of this section, we
give a step-by-step argument as follows. By the condition (d) of
Definition 3.1, Proposition1.1 entails that $A$ has the PBW
$K$-basis, i.e. the condition (1) of Theorem 2.6 is satisfied. By
the definition of ${\cal R}$,  $A$ satisfies the condition (2) of
Theorem 2.6. Since every binomial skew polynomial ring is a domain
as established in [G-I1,2],  the condition (3) of Theorem 2.6 is
satisfied by $A$ as well. }\QED\v5

Let $A=\KX /I=K[a_1,\ldots ,a_n]$ be a binomial skew polynomial ring
as defined above, where $a_i$ stands for the coset $X_i+I$ in $\KX
/I$, $1\le i\le n$. Then $A$ is clearly a {\it skew 2-nomial
algebra} in the sense of [Li3], and it follows from the conditions
(a) -- (d) of  Definition 3.1 that numerous  $A$ may fall into the
case of  [Li3, Theorem 2.2]. If it is the case, then $A$ has the
lexicographic ordering $\prec_{lex}$ on the PBW $K$-basis $\B$ such
that
$a_n\prec_{lex}a_{n-1}\prec_{lex}\cdots\prec_{lex}a_2\prec_{lex}a_1$,
and $A$ has a left Gr\"obner basis theory with respect to
$\prec_{lex}$, i.e., every left ideal $L$ of $A$ has a finite left
Gr\"obner basis in the sense that if $f\in L$ and $f\ne 0$, then
there is a $g\in\G$ such that $\LM (g)|\LM (f)$ (a left division),
where $\LM(~)$ stands for taking the leading monomial of elements in
$A$ with respect to  $\prec_{lex}$ (note that a binomial skew
polynomial ring $A$ is Noetherian and hence a left Gr\"obner basis
of $L$ is finite). Thus, since for each $s=1,2,\ldots ,n-1$  the
monomial ordering $\prec_{lex}$ is an elimination ordering of type
$s$ with respect to ${\bf V}(U_s)$, i.e., $f\in A$ with the leading
monomial $\LM (f)\in {\bf V}(U_s)$ implies $f\in {\bf V}(U_s)$,
where $U_s=\{ a_{n},\ldots ,a_{n-s+1}\}$, ${\bf
V}(U_s)=K\hbox{-span}\{ a_{n-s+1}^{\alpha_1}\cdots a_{n}^{\alpha_s}~
|~ (\alpha_1,...,\alpha_s)\in\NZ^s\}$, if $L$ is a nonzero left
ideal of $A$, then GK.dim$(A/L)=d<$ GK.dim$A=n$, ${\bf
V}(U_{d+1})\cap L\ne \{0\}$ by Theorem 3.2 and Lemma 2.1, and a left
Gr\"obner basis $\G$ of $L$ with respect to $\prec_{lex}$ will
contain a nonzero element of ${\bf V}(U_{d+1})$. Therefore,
realizing the elimination property of $L$ (as described in Lemma
2.1) may become possible in case computing a left Gr\"obner basis
for $L$ under $\prec_{lex}$ is algorithmically feasible.\v5

By Definition 3.1 it is also clear that a binomial skew polynomial
ring $A$ is not necessarily a quadric solvable polynomial algebra in
the sense of [Li1, CH.III, Section 2]. Hence, with Theorem 3.2 we
end this section by concluding that the class of algebras satisfying
the three conditions of Theorem 2.6  properly contains the class of
all quadric solvable polynomial algebras in the sense of [Li1,
CH.III, Section 2]. This illustrates Remark (ii) given at the end of
Section 2. \v5

\section*{4. An Application to Solvable Polynomial Algebras}
As we have seen  from the introduction section and  Section 2, that
Elimination Lemma 2.1 holds true for every  nonzero left ideal of a
quadric solvable polynomial algebra in the sense of [Li1, CH.III,
Section 2]. In this section, we show that, indeed, Elimination Lemma
2.1 holds true (in the sense of Definition 2.2) for every nonzero
left ideal of an arbitrary solvable polynomial algebra in the sense
of [K-RW].  To this end, we first recall the
following{\parindent=0pt\v5

{\bf 4.1. Definition} ([K-RW], [LW]) Let $A=K[a_1,\ldots ,a_n]$ be a
finitely generated $K$-algebra. Suppose that  $A$ has the PBW
$K$-basis $\B=\{a^{\alpha}=a_1^{\alpha_1}\cdots
a_n^{\alpha_n}~|~\alpha =(\alpha_1,\ldots ,\alpha_n)\in\NZ^n\}$,
and that $\prec$ is a (two-sided) monomial ordering on $\B$. If for
all $a^{\alpha}=a_1^{\alpha_1}\cdots a_n^{\alpha_n}$,
$a^{\beta}=a_1^{\beta_1}\cdots a^{\beta_n}_n\in\B$, the following
holds:
$$\begin{array}{rcl} a^{\alpha}a^{\beta}&=&\lambda_{\alpha ,\beta}a^{\alpha +\beta}+f_{\alpha ,\beta},\\
&{~}&\hbox{where}~\lambda_{\alpha ,\beta}\in K^*,~a^{\alpha
+\beta}=a_1^{\alpha_1+\beta_1}\cdots a_n^{\alpha_n+\beta_n},~\hbox{and either}~f_{\alpha ,\beta}=0~\hbox{or}\\
&{~}&f_{\alpha ,\beta}=\sum_q\mu_qa^{\gamma (q)}\in
A~\hbox{with}~\mu_q\in K,~a^{\gamma (q)}\in\B,~\hbox{satisfying}~\LM
(f_{\alpha ,\beta})\prec a^{\alpha +\beta},\end{array}$$ where $\LM
(f_{\alpha ,\beta})$ stands for the leading monomial of $f_{\alpha
,\beta}$ with respect to $\prec$, then $A$ is called a {\it solvable
polynomial algebra}.\v5

{\bf 4.2. Proposition}  [Li5, Theorem 2.1] Let $A=K[a_1,\ldots
,a_n]$ be a finitely generated $K$-algebra, and let $\KS =K\langle
X_1,\ldots ,X_n\rangle$ be the free $K$-algebras with the standard
$K$-basis $\mathbb{B}=\{ 1\}\cup\{X_{i_1}\cdots X_{i_s}~|~X_{i_j}\in
X,~s\ge 1\}$. With notations as before, the following two statements
are equivalent:\par

(i) $A$ is a solvable polynomial algebra in the sense of Definition
4.1.\par

(ii) $A\cong \OV A=\KS /I$ via the $K$-algebra epimorphism $\pi_1$:
$\KS \r A$ with $\pi_1(X_i)=a_i$, $1\le i\le n$, $I=$ Ker$\pi_1$,
satisfying  {\parindent=1.3truecm

\item{(a)} with respect to some monomial ordering $\prec_{_X}$ on $\mathbb{B}$, the ideal $I$ has a
finite Gr\"obner basis $G$ and the reduced Gr\"obner basis of $I$ is
of the form
$$\G =\left\{ g_{ji}=X_jX_i-\lambda_{ji}X_iX_j-F_{ji}
~\left |~\begin{array}{l} \LM (g_{ji})=X_jX_i,\\ 1\le i<j\le
n\end{array}\right. \right\} $$ where $\lambda_{ji}\in K^*$,
$\mu^{ji}_q\in K$, and
$F_{ji}=\sum_q\mu^{ji}_qX_1^{\alpha_{1q}}X_2^{\alpha_{2q}}\cdots
X_n^{\alpha_{nq}}$ with $(\alpha_{1q},\alpha_{2q},\ldots
,\alpha_{nq})\in\NZ^n$, thereby $\B =\{ \OV X_1^{\alpha_1}\OV
X_2^{\alpha_2}\cdots \OV X_n^{\alpha_n}~|~$ $\alpha_j\in\NZ\}$ forms
a PBW $K$-basis for $\OV A$, where each $\OV X_i$ denotes the coset
of $I$ represented by $X_i$ in $\OV A$; and

\item{(b)} there is a (two-sided) monomial ordering
$\prec$ on $\B$ such that $\LM (\OV{F}_{ji})\prec \OV X_i\OV X_j$
whenever $\OV F_{ji}\ne 0$, where $\OV F_{ji}=\sum_q\mu^{ji}_q\OV
X_1^{\alpha_{1i}}\OV X_2^{\alpha_{2i}}\cdots \OV X_n^{\alpha_{ni}}$,
$1\le i<j\le n$. \par\QED}}\v5

In conclusion, we derive the next{\parindent=0pt\v5

{\bf 4.3. Theorem}  Let $A=K[a_1,\ldots ,a_n]$ be any solvable
polynomial algebra in the sense of Definition 4.1. Then  Elimination
Lemma 2.1 holds true (in the sense of Definition 2.2) for every
nonzero left ideal $L$ of $A$.\vskip 6pt

{\bf Proof} First note that every solvable polynomial algebra $A$
has the PBW $K$-basis $\B$ by Definition 4.1 (or by Proposition 1.1
and Theorem 4.2). Moreover, it follows from  Proposition 1.2 and
Proposition 4.2 that  GK.dim$A=n$ ($=$ the number of generators of
$A$). As also we know that $A$ is a domain by [K-RW]. Hence Lemma
2.3   entails that GK.dim$(A/L)<$ GK.dim$A=n$ holds for every
nonzero left ideal of $A$. Therefore, we conclude that Elimination
Lemma 2.1 holds true (in the sense of Definition 2.2) for every
nonzero left ideal $L$ of $A$.  \QED}\v5

Comparing with  Theorem 2.6, we see that Theorem 4.3 did not require
$A$ satisfies the condition (2) of Theorem 2.6 (or equivalently,
Theorem 4.3 did not require $A$ is an $\NZ$-filtered algebra with
respect to the filtration $\F A$ determined by the PBW $K$-basis
$\B$). Indeed, the example presented below illustrates that a
solvable polynomial algebra $A$ may not necessarily satisfy the
condition (2) of Theorem 2.6 (or equivalently, $A$ may not be  an
$\NZ$-filtered algebra with respect to the filtration $\F A$
determined by the PBW $K$-basis $\B$). {\parindent=0pt\v5

{\bf Example}  [Li5, Example 1]  Considering the positive-degree
function $d$ on the free $K$-algebra $\KS =K\langle X_1,X_2,
X_3\rangle$ such that $d(X_1)=2$, $d(X_2)=1$, and $d(X_3)=4$, let
$I$ be the ideal of $\KS$ generated by the subset $\G$ consisiting
of
$$\begin{array}{l} g_1=X_1X_2- X_2X_1,\\
g_2=X_3X_1-\lambda X_1X_3-\mu X_3X_2^2-f(X_2),\\
g_3=X_3X_2- X_2X_3,\end{array}$$ where $\lambda\in K^*$, $\mu\in K$,
$f(X_2)$ is a polynomial in $X_2$ with $d(f(X_2))\le 6$, or
$f(X_2)=0$. The following statements hold.   \par

(1) $\G$ forms a Gr\"obner basis for $I$  with respect to the graded
lexicographic ordering $X_2\prec_{grlex}X_1\prec_{grlex}X_3$,  such
that the three generators of $I$ have the leading monomials $\LM
(g_1)=X_1X_2$, $\LM (g_2)=X_3X_1$, and $\LM (g_3)=X_3X_2$. \par

(2) With respect to the fixed $\prec_{grlex}$ in (1), the reduced
Gr\"obner basis $\G '$ of $I$ consists of
$$\begin{array}{l} g_1=X_1X_2- X_2X_1,\\
g_2=X_3X_1- \lambda X_1X_3-\mu X_2^2X_3-f(X_2),\\
g_3=X_3X_2- X_2X_3,\end{array}$$\par

(3) Writing $A=K[a_1,a_2, a_3]$ for the quotient algebra $\KS /I$,
where $a_1$, $a_2$ and $a_3$ denote the cosets $X_1+I$, $X_2+I$ and
$X_3+I$ in $\KS /I$ respectively, it follows that $A$ has the PBW
basis $\B =\{
a^{\alpha}=a_2^{\alpha_2}a_1^{\alpha_1}a_3^{\alpha_3}~|~\alpha
=(\alpha_2,\alpha_1,\alpha_3)\in\NZ^3\}$. Noticing that $
a_2a_1=a_1a_2$, it is clear that $\B ' =\{
a^{\alpha}=a_1^{\alpha_1}a_2^{\alpha_2}a_3^{\alpha_3}~|~\alpha
=(\alpha_1,\alpha_2,\alpha_3)\in\NZ^3\}$ is also a PBW basis for
$A$. Since $a_3a_1=\lambda a_1a_3+\mu a_2^2a_3+f(a_2)$, where
$f(a_2)\in K$-span$\{ 1,a_2,a_2^2,\ldots ,a_2^6\}$, we see that $A$
has the monomial ordering $\prec_{lex}$ on $\B '$ such that
$a_3\prec_{lex}a_2\prec_{lex}a_1$ and $\LM (\mu
a_2^2a_3+f(a_2))\prec_{lex}a_1a_3$, thereby $A$ is turned into a
solvable polynomial algebra with respect to $\prec_{lex}$. Also if
we use the positive-degree function $d$ on $\B '$ such that
$d(a_1)=2$, $d(a_2)=1$, and $d(a_3)=4$, then $A$ has another
monomial ordering on $\B '$, namely the graded lexicographic
ordering $\prec_{grlex}$ such that
$a_3\prec_{grlex}a_2\prec_{grlex}a_1$ and $\LM (\mu
a_2^2a_3+f(a_2))\prec_{grlex}a_1a_3$, thereby $A$ is turned into a
solvable polynomial algebra with respect to $\prec_{grlex}$.\par

(4) Consider the $K$-algebra $A=K[a_1,a_2, a_3]$ as presented in the
above (3), and consider the positive-degree function $d$ on the PBW
$K$-basis $\B '$ of $A$ such that $d(a_i)=1$, $1\le i\le 3$. Then
$A$ is not an $\NZ$-filtered algebra with respect to the filtration
$\F A$ determined by  $\B '$, because $a_3a_1=\lambda a_1a_3+\mu
a_2^2a_3+f(a_2)$ implies $\F_1A\F_1A\not\subset \F_2 A$. Therefore,
the solvable polynomial algebra $A$ does not satisfies the condition
(2) of Theorem 2.6 (see also Lemma 2.4), illustrating Remark (iii)
given at the end of Section2.} \v5

Finally, note that a noncommutative Gr\"obner basis theory works
effectively for every solvable polynomial algebra $A=K[a_1,\ldots
,a_n]$ in the sense of [K-RW], that is,  a noncommutative
Buchberger's algorithm works very well in the sense that if a finite
generating set of a left ideal $L$ is given (note that $A$ is
Noetherian), then running the noncommutative Buchberger's algorithm
with respect to a monomial ordering $\prec$ will produce a finite
left Gr\"obner basis $\G$ for $L$. Thus,  since GK.dim$A/L=d<n$ by
Theorem 4.3, if $U=\{ a_{i_1},...,a_{i_{d+1}}\}$ with
$i_1<i_2<\cdots <i_{d+1}$, and if $A$ admits an elimination ordering
$\prec_{d+1}$ of type $d+1$ with respect to ${\bf
V}(U)=K$-span$\{a_{i_1}^{\alpha_1}\cdots
a_{i_{d+1}}^{\alpha_{d+1}}~|~(\alpha_1,\ldots
,\alpha_{d+1})\in\NZ^{d+1}\}$, i.e., $f\in A$ with the leading
monomial $\LM (f)\in {\bf V}(U)$ implies $f\in {\bf V}(U)$, then the
noncommutative Buchberger's algorithm will produce a left Gr\"obner
~basis $\G$ of $L$ such that $\G$ contains a nonzero element of
${\bf V}(U)$. This shows that the elimination property of $L$ (as
described in Lemma 2.1) may be realized via computing  a left
Gr\"obner basis of $L$ with respect to a suitable elimination
ordering.\v5

\centerline{\bf References} {\parindent=.8truecm\vskip 8pt

\item{[BW]} T.~Becker and V.~Weispfenning, {\it Gr\"obner Bases},
Springer-Verlag, 1993.

\item{[Ch]} F.~Chyzak, Holonomic systems and automatic proving of
identities, {\it Research Report} 2371,  Institute National de
Recherche en Informatique et en Automatique, 1994.

\item{[CS]} F. Chyzak and B. Salvy, Noncommutative elimination in Ore
algebras proves multivariate identities, {\it J. Symbolic Comput}.,
26(1998), 187--227.

\item{[G-I1]} T. Gateva-Ivanova, Skew polynomial rings with binomial
relations. {\it J. Algebra},  185(1996), 710--753.

\item{[G-I2]} T. Gateva-Ivanova, Binomial skew polynomial rings, Artin-Schelter regularity, and binomial
solutions of the Yang-Baxter equation.  {Serdica Math. J}. 30(2004),
431--470

\item{[G-I3]} T. Gateva-Ivanova, Quadratic algebras, Yang-Baxter equation, and Artin-Schelter
regularity£¬{\it Advances in Mathematics}, 230(2012), 2152¨C2175.

\item{[Gr]} E. L. Green,  Noncommutative
Gr\"obner bases and projective resolutions.  In: {\it Computational
Methods for Representations of Groups and Algebras} (Michler,
Schneider, eds), Proceedings of the Euroconference, Essen, 1997.
Progress in Mathematics, Vol. 173, Basel, Birkha¡§user Verlag, 1999,
29--60.\par

\item{[Gr\"o]} W.~Gr\"obner, {\it Algebraic Geometrie} I, II,
Bibliographisches Institut, Mannheim, 1968, 1970.

\item{[Hu]} J. E. Humphreys, {\it Introduction to Lie Algebras and Representation Theory}. Springer, 1972.

\item{[KL]} G.R. Krause and T.H. Lenagan, {\it Growth of Algebras and
Gelfand-Kirillov Dimension}. Graduate Studies in Mathematics.
American Mathematical Society, 1991.

\item{[KW]} H. Kredel and V.~Weispfenning, Computing dimension and
independent sets for polynomial ideals. In: {\it Computational
Aspects of Commutative Algebra}, from a special issue of the Journal
of Symbolic Computation (L.~Robbiano ed), Academic Press, 1989,
97--113.

\item{[KR-W]} A. Kandri-Rody, V.~Weispfenning, Non-commutative
Gr\"obner bases in algebras of solvable type.  {\it J. Symbolic
Comput.}, 9(1990), 1--26.\par

\item{[Li1]} H. Li, {\it Noncommutative Gr\"obner Bases and
Filtered-Graded Transfer}. Lecture Notes in Mathematics, Vol. 1795,
Springer, 2002.\par

\item{[Li2]} H. Li, $\Gamma$-leading homogeneous algebras and Gr\"obner
bases. In: {\it Recent Developments in Algebra and Related Areas}
(F. Li and C. Dong eds.), Advanced Lectures in Mathematics, Vol. 8,
International Press \& Higher Education Press, Boston-Beijing, 2009,
155--200. (This is a  strengthened version of
\textsf{arXiv:math.RA/0609583}, http://arXiv.org)

\item{[Li3]} H. Li, Looking for Gr\"obner basis theory for (almost) skew 2-nomial
algebras. {\it J. Symbolic Comput.}, 45(2010), 918¨C942. Available 
at \textsf{arXiv:0808.1477 [math.RA]}, http://arXiv.org

\item{[Li4]} H. Li, {\it Gr\"obner Bases in Ring Theory}. World Scientific Publishing Company, 2011.\par

\item{[Li5]} H. Li, A Note on solvable polynomial algebras. {\it Computer Science Journal of Moldova}, vol.22, 1(64), 2014, 99 -- 109.
Available at  \textsf{arXiv:1212.5988 [math.RA]}, http://arXiv.org

\item{[LW]} H. Li.  Y. Wu, Filtered-graded transfer of
Gr\"obner basis computation in solvable polynomial algebras. {\it
Communications in Algebra}, 28(1), 2000, 15--32.\par

\item{[MR]} J.C.~McConnell and J.C.~Robson, {\it Noncommutative
Noetherian Rings}, John Wiley \& Sons, 1987.

\item{[Mor]} T. Mora, An introduction to commutative and noncommutative
Gr\"obner Bases. {\it Theoretic Computer Science}, 134, 1994,
131--173.

\item{[PWZ]} M. Petkov¡¦sek, H. Wilf and D. Zeilberger, {\it $A = B  $}. A.K.
Peters, Ltd. 1996.

\item{[SW]} A.V. Shepler and S. Witherspoon,
Poincar\'e-Birkhoff-Witt Theorems. In: {\it Commutative Algebra and
Noncommutative Algebraic Geometry} (D. Eisenbud, S.B. Iyengar, A.K.
Singh, J.T. Stafford, and M. Van den Bergh, eds), Mathematical
Sciences Research Institute Proceedings, Vol. 1, Cambridge Univ.
Press, 2015.

\item{[Uf]} V. Ufnarovski, A growth criterion for graphs and algebras
defined by words. {\it Mat. Zametki}, 31(1982), 465--472 (in
Russian); English translation: {\it Math. Notes}, 37(1982),
238--241.

\item{[WZ]} H.S. Wilf and D. Zeilberger, An algorithmic proof theory
for hypergeometric (ordinary and ¡°q¡±) multisum/integral 
identities, {\it Invent. Math}., 108(1992), 575 -- 633.

\item{[Zei1]} D. Zeilberger, A holonomic system approach to special
function identities, {\it J. Comput. Appl. Math}., 32(1990), 321 --
368.

\item{[Zei2]} D. Zeilberger, The method of creative telescoping, {\it J.
Symbolic Comput}., 11(1991), 195¨C204.\par}

\end{document}